\documentclass[12pt]{amsart}
\usepackage{amsmath, amssymb, latexsym, amsthm, mathrsfs}
\usepackage[all]{xy}

\newcommand{\ed}{

\end{document}}

\def\su{\subseteq}
\def\om{\omega}

\def\forces{{\;\Vdash}}
\newcommand\name[1]{\mathring{#1}}
\newcommand{\Pa}[8]{\bibitem{#1} {#2}, \emph{#3}, {#4} \textbf{#5} ({#6}), {#7}--{#8}.}

\newcommand{\Setting}[7]{\xymatrix@R=4pt@C=7pt{#1\ar@{-}[r]&#2\ar@{-}[r]&#3\\&#4\ar@{-}[u]\\
#5\ar@{-}[uu]\ar@{-}[r] & #6\ar@{-}[u]\ar@{-}[r] & #7\ar@{-}[uu]}}

\newcommand{\bbP}{\mathbb{P}}

\def\rmand{\mbox{ and }}

\newcommand{\bq}{\begin{quote}}
\newcommand{\eq}{\end{quote}}

\newcommand{\N}{\w} 
\newcommand{\NN}{{\N^{\N}}}

\newcommand{\roth}{{[\w]^{\w}}}
\newcommand{\Fin}{{[\w]^{<\w}}}

\newcommand{\sseq}[1]{\{#1 : n\in\N\}}

\newcommand{\op}{\operatorname}

\newcommand{\cB}{\mathcal{B}}
\newcommand{\cC}{\mathcal{C}}
\newcommand{\cH}{\mathcal{H}}

\newcommand{\cA}{\mathcal{A}}

\newcommand{\cD}{\mathcal{D}}
\newcommand{\cF}{\mathcal{F}}

\newcommand{\cM}{\mathcal{M}}
\newcommand{\cN}{\mathcal{N}}

\newcommand{\bbQ}{\mathbb{Q}}
\newcommand{\bbR}{\mathbb{R}}
\newcommand{\cU}{\mathcal{U}}
\newcommand{\Union}{\bigcup}
\newcommand{\cV}{\mathcal{V}}
\newcommand{\cW}{\mathcal{W}}

\newcommand{\Impl}{\Rightarrow}
\long\def\forget#1\forgotten{}
\newcommand{\ft}{\mathfrak{t}}
\newcommand{\fb}{\mathfrak{b}}
\newcommand{\fc}{\mathfrak{c}}
\newcommand{\fd}{\mathfrak{d}}

\newcommand{\oo}{\infty}

\newcommand{\fh}{\mathfrak{h}}

\newcommand{\w}{\omega}
\newcommand{\x}{\times}

\newcommand{\nin}{\notin}
\def\res{\upharpoonright}
\newcommand\concat{{\hat{\phantom{aa}}}}

\newcommand{\sbst}{\subseteq}

\newcommand{\sm}{\setminus}
\newcommand{\as}{\subseteq^*}

\newcommand{\cov}{\op{cov}}

\newcommand{\cof}{\op{cof}}

\newcommand{\non}{\op{non}}

\newtheorem{thm}{Theorem}[section]
\newcommand{\bthm}{\begin{thm}} \newcommand{\ethm}{\end{thm}}
\newtheorem{prop}[thm]{Proposition}
\newcommand{\bprp}{\begin{prop}} \newcommand{\eprp}{\end{prop}}
\newtheorem{fact}[thm]{Fact}
\newcommand{\bfct}{\begin{fact}} \newcommand{\efct}{\end{fact}}
\newtheorem{prob}[thm]{Problem}
\newcommand{\bprb}{\begin{prob}} \newcommand{\eprb}{\end{prob}}
\newtheorem{lem}[thm]{Lemma}
\newcommand{\blem}{\begin{lem}} \newcommand{\elem}{\end{lem}}
\newtheorem{claim}[thm]{Claim}
\newcommand{\bclm}{\begin{claim}} \newcommand{\eclm}{\end{claim}}
\newtheorem{cor}[thm]{Corollary}
\newcommand{\bcor}{\begin{cor}} \newcommand{\ecor}{\end{cor}}
\newtheorem{conj}[thm]{Conjecture}
\newcommand{\bcnj}{\begin{conj}} \newcommand{\ecnj}{\end{conj}}
\theoremstyle{definition}
\newtheorem{defn}[thm]{Definition}
\newcommand{\bdfn}{\begin{defn}} \newcommand{\edfn}{\end{defn}}
\theoremstyle{remark}
\newtheorem{rem}[thm]{Remark}
\newcommand{\brem}{\begin{rem}} \newcommand{\erem}{\end{rem}}
\newtheorem{cnv}[thm]{Convention}
\newcommand{\bcnv}{\begin{cnv}} \newcommand{\ecnv}{\end{cnv}}
\newtheorem{exam}[thm]{Example}
\newcommand{\bexm}{\begin{exam}} \newcommand{\eexm}{\end{exam}}
\newcommand{\bpf}{\begin{proof}} \newcommand{\epf}{\end{proof}}
\newcommand{\be}{\begin{enumerate}}
\newcommand{\ee}{\end{enumerate}}
\newcommand{\bi}{\begin{itemize}}
\newcommand{\itm}{\item}
\newcommand{\ei}{\end{itemize}}

\forget
\forgotten

\newcommand{\sone}{\mathsf{S}_1}

\newcommand{\scof}{\mathsf{S}_\mathrm{cof}}
\newcommand{\ufin}{\mathsf{U}_\mathrm{fin}}

\title{Point-cofinite covers in the Laver model}
\author[Arnold Miller]{Arnold W. Miller}
\address[Miller]{Department of Mathematics, University of Wisconsin-Madison, Van Vleck Hall 480 Lincoln Drive, Madison, Wisconsin 53706-1388, USA}
\email{miller@math.wisc.edu}
\urladdr{http://www.math.wisc.edu/\~{}miller/}
\author{Boaz Tsaban}
\address[Tsaban]{Department of Mathematics, Bar-Ilan University, Ramat-Gan 52900, Israel}
\email{tsaban@math.biu.ac.il}
\urladdr{http://www.cs.biu.ac.il/\~{}tsaban}

\begin{document}

\newcommand{\sgg}{\sone(\Gamma,\Gamma)}

\begin{abstract}
Let $\sgg$ be the statement: For each sequence of point-cofinite open covers,
one can pick one element from each cover and obtain a point-cofinite cover.
$\fb$ is the minimal cardinality of a set of reals not satisfying $\sgg$.
We prove the following assertions:
\be
\itm If there is an unbounded tower, then there are sets of reals of cardinality $\fb$, satisfying $\sgg$.
\itm It is consistent that all sets of reals satisfying $\sgg$ have cardinality smaller than $\fb$.
\ee
These results can also be formulated as dealing with Arhangel'ski\u{\i}'s property $\alpha_2$
for spaces of continuous real-valued functions.

The main technical result is that in Laver's model,
each set of reals of cardinality $\fb$ has an unbounded
Borel image in the Baire space $\NN$.
\end{abstract}

\maketitle

\section{Background}

Let $P$ be a nontrivial property of sets of reals.
The \emph{critical cardinality} of $P$,
denoted $\non(P)$, is the minimal cardinality of a set of reals not satisfying $P$.
A natural question is whether there is a set of reals of cardinality at least $\non(P)$,
which satisfies $P$, i.e., a \emph{nontrivial} example.


We consider the following property.
Let $X$ be a set of reals. $\cU$ is a \emph{point-cofinite}
cover of $X$ if $\cU$ is infinite,
and for each $x\in X$, $\{U\in\cU : x\in U\}$ is a cofinite subset of $\cU$.\footnote{Historically, point-cofinite covers were named
\emph{$\gamma$-covers}, since they are related to a property numbered $\gamma$ in a list from $\alpha$ to $\epsilon$
in the seminal paper \cite{GN} of Gerlits and Nagy.}
Having $X$ fixed in the background, let $\Gamma$ be the family of all point-cofinite \emph{open} covers of $X$.
The following properties were introduced by Hurewicz \cite{Hure25}, Tsaban \cite{MHP}, and Scheepers \cite{coc1}, respectively.
\begin{description}
\item[$\ufin(\Gamma,\Gamma)$] For all $\cU_0,\cU_1,\dots\in\Gamma$, none containing
a finite subcover, there are finite
$\cF_0\sbst\cU_0,\cF_1\sbst\cU_1,\dots$ such that
$\sseq{\Union\cF_n}\in\Gamma$.
\item[$\mathsf{U}_2(\Gamma,\Gamma)$] For all $\cU_0,\cU_1,\dots\in\Gamma$, there are
$\cF_0\sbst\cU_0,\cF_1\sbst\cU_1,\dots$ such that $|\cF_n|=2$ for all $n$, and
$\sseq{\Union\cF_n}\in\Gamma$.
\item[$\sgg$] For all $\cU_0,\cU_1,\dots\in\Gamma$, there are
$U_0\in\cU_0,U_1\in\cU_1,\dots$ such that $\sseq{U_n}\in\Gamma$.
\end{description}
Clearly, $\sgg$ implies $\mathsf{U}_2(\Gamma,\Gamma)$, which in turn implies $\ufin(\Gamma,\Gamma)$.
None of these implications is reversible in ZFC \cite{MHP}.
The critical cardinality of all three properties is $\fb$ \cite{coc2}.\footnote{Blass's survey \cite{BlassHBK}
is a good reference for the definitions and details about the special cardinals mentioned in this paper.}

Bartoszy\'nski and Shelah \cite{BaSh01} proved that
there are, provably in ZFC, totally imperfect sets of reals of cardinality $\fb$
satisfying the Hurewicz property $\ufin(\Gamma,\Gamma)$. Tsaban proved the same assertion for $\mathsf{U}_2(\Gamma,\Gamma)$ \cite{MHP}.
These sets satisfy $\ufin(\Gamma,\Gamma)$ in all finite powers \cite{ideals}.

We show that in order to obtain similar results for $\sgg$, hypotheses beyond ZFC are necessary.

\section{Constructions}

We show that certain weak (but not provable in ZFC) hypotheses suffice to have nontrivial $\sgg$ sets,
even ones which possess this property in all finite powers.

\bdfn
A \emph{tower} of cardinality $\kappa$ is a set $T\sbst\roth$ which
can be enumerated
bijectively as $\{x_\alpha : \alpha<\kappa\}$,
such that for all $\alpha<\beta<\kappa$, $x_\beta\as x_\alpha$.

A set $T\sbst\roth$ is \emph{unbounded} if the
set of its enumeration functions are unbounded, i.e., for any $g\in\om^\om$
there is an $x\in T$ such that for infinitely many $n$,
$g(n)$ is less than the $n$-th element of $x$.
\edfn

Scheepers \cite{alpha_i} proved that if $\ft=\fb$, then there is a set of reals
of cardinality $\fb$, satisfying $\sgg$.
If $\ft=\fb$, then there is an unbounded tower of cardinality $\fb$,
but the latter assumption is weaker.

\blem[folklore]\label{bltd}
If $\fb<\fd$, then there is an unbounded tower of cardinality $\fb$.
\elem
\bpf
Let $B=\{b_\alpha : \alpha<\fb\}\sbst\NN$ be a $\fb$-scale, that is,
each $b_\alpha$ is increasing, $b_\alpha\le^* b_\beta$
for all $\alpha<\beta<\fb$,
and $B$ is unbounded.

As $|B|<\fd$, $B$ is not dominating. Let $g\in\NN$ exemplify that.
For each $\alpha<\fb$, let $x_\alpha=\{n : b_\alpha(n)\le g(n)\}$.
Then $T=\{x_\alpha : \alpha<\fb\}$ is an unbounded tower: Clearly, $x_\beta\as x_\alpha$
for $\alpha<\beta$. Assume that $T$ is bounded, and let $f\in\NN$ exemplify that.
For each $\alpha$, writing $x_\alpha(n)$ for the $n$-th element of $x_\alpha$:
$$b_\alpha(n)\le b_\alpha(x_\alpha(n))\le g(x_\alpha(n))\le g(f(n))$$
for all but finitely many $n$. Thus, $g\circ f$ shows that $B$ is bounded.
A contradiction.
\epf

\bthm\label{main}
If there is an unbounded tower (of any cardinality), then
there is a set of reals $X$
of cardinality $\fb$, which satisfies $\sgg$.
\ethm

Theorem \ref{main} follows from the following two propositions.

\bprp\label{anytower}
If there is an unbounded tower, then there is one of cardinality $\fb$.
\eprp
\bpf
By Lemma \ref{bltd}, it remains to consider the case $\fb=\fd$.
Let $T$ be an unbounded tower of cardinality $\kappa$.
Let $\{f_\alpha : \alpha<\fb\}\sbst\NN$ be dominating.
For each $\alpha<\fb$, pick $x_\alpha\in T$ which is not bounded by $f_\alpha$.
$\{x_\alpha : \alpha<\fb\}$ is unbounded, being unbounded in a dominating family.
\epf

Define a topology on $P(\om)$ by identifying $P(\om)$ with the Cantor space $2^\om$,
via characteristic functions.
Scheepers's mentioned proof actually establishes the following result, to which
we give an alternative proof.

\bprp[essentially, Scheepers \cite{alpha_i}]\label{Sch}
For each unbounded tower $T$ of cardinality $\fb$,
$T\cup\Fin$ satisfies $\sgg$.
\eprp
\bpf
Let $T=\{x_\alpha : \alpha<\fb\}$ be an unbounded tower of cardinality $\fb$.
For each $\alpha$, let $X_\alpha=\{x_\beta : \beta<\alpha\}\cup\Fin$.
Let $\cU_0,\cU_1,\dots$ be point-cofinite open covers of $X_\fb=T\cup\Fin$.
We may assume that each $\cU_n$ is countable and that
$\cU_i\cap\cU_j=\emptyset$ whenever $i\neq j$.

By the proof of Lemma 1.2 of \cite{GM}, for each $k$ there are distinct $U^k_0,U^k_1,\dots\in\cU_k$,
and an increasing sequence $m^k_0<m^k_1<\dots$, such that for each $n$ and $k$,
$$\{ x\su\om \;:\; x\cap (m^k_n,m^k_{n+1})=\emptyset\}\su U^k_n.$$
As $T$ is unbounded, there is $\alpha<\fb$ such that
for each $k$, $I_k = \{n : x_\alpha\cap (m^k_n,m^k_{n+1})=\emptyset\}$
is infinite.

For each $k$, $\{U^k_n : n\in \om\}$ is an infinite subset
of $\cU_k$, and thus
a point-cofinite cover of $X_\alpha$.
As $|X_\alpha|<\fb$, there is $f\in\NN$ such that
$$\forall x\in X_\alpha\;\exists k_0\;\forall k\geq k_0\;
\forall n>f(k)\;\; x\in U^k_n.$$
For each $k$, pick $n_k\in I_k$ such that $n_k>f(k)$,

We claim that $\{U^k_{n_k} : k\in\N\}$ is a point-cofinite cover
of $X_\fb$:
If $x\in X_\alpha$, then $x\in U^k_{n_k}$ for all but finitely many $k$, because
$n_k> f(k)$ for all $k$.
If $x=x_\beta$, $\beta\ge\alpha$, then $x\as x_\alpha$.
For each large enough $k$, $m^k_{n_k}$ is large enough, so that
$x\cap (m^k_{n_k},m^k_{n_k+1})\sbst x_\alpha\cap (m^k_{n_k},m^k_{n_k+1})
=\emptyset$,
and thus $x\in U^k_{n_k}$.
\epf

\brem\label{pst}
Zdomskyy points out that for the proof to go through, it suffices that $\{x_\alpha : \alpha<\fb\}$ is such that
there is an unbounded $\{y_\alpha : \alpha<\fb\}\sbst\roth$ such that for each $\alpha$,
$x_\alpha$ is a pseudointersection of $\{y_\beta : \beta<\alpha\}$.
We do not know whether the assertion mentioned here is weaker than the existence of an unbounded tower.
\erem

We now turn to nontrivial examples of sets satisfying $\sgg$ in all finite powers.
In general, $\sgg$ is not preserved by taking finite powers \cite{coc2}, and
we use a slightly stronger hypothesis in our construction.

\bdfn
Let $\fb_0$ be the additivity number of $\sgg$, that is, the minimum cardinality
of a family $\cF$ of sets of reals, each satisfying $\sgg$, such that the union of all members
of $\cF$ does not satisfy $\sgg$.
\edfn

$\ft\le\fh$, and Scheepers proved that $\fh\le\fb_0\le\fb$ \cite{wqn}. It follows from
Theorem \ref{mainthm} that consistently, $\fh<\fb_0=\fb$. It is open whether $\fb_0=\fb$ is provable.
If $\ft=\fb$ or $\fh=\fb<\fd$, then there is an unbounded tower of cardinality $\fb_0$.

\bthm
For each unbounded tower $T$ of cardinality $\fb_0$,
all finite powers of $T\cup\Fin$ satisfy $\sgg$.
\ethm
\bpf
We say that $\cU$ is an \emph{$\w$-cover} of $X$ if no member of $\cU$ contains $X$
as a subset, but each finite subset of $X$ is contained in some member of $\cU$.
We need a multidimensional version of Lemma 1.2 of \cite{GM}.

\blem\label{GM++}
Assume that $\Fin\sbst X\sbst P(\N)$, and let $e\in\N$.
For each open $\w$-cover $\cU$ of $X^e$, there are $m_0<m_1<\dots$ and
 $U_0,U_1,\dots\in\cU$,
such that for all $x_0,\dots,x_{e-1}\sbst\N$, $(x_0,\dots,x_{e-1})\in U_n$ whenever
$x_i\cap (m_n,m_{n+1})=\emptyset$ for all $i<e$.
\elem
\bpf
As $\cU$ is an open $\w$-cover of $X^e$, there is an open $\w$-cover $\cV$ of $X$ such that $\{V^e : V\in\cV\}$ refines $\cU$ \cite{coc2}.

Let $m_0=0$. For each $n\ge 0$:
Assume that $V_0,\dots,V_{n-1}\in\cV$ are given, and $U_0,\dots,U_{n-1}\in\cU$
are such that $V_i^e\sbst U_i$ for all $i<n$.
Fix a finite $F\sbst X$ such that $F^e$ is not contained
in any of the sets $U_0,\dots,U_{n-1}$.
As $\cV$ is an $\w$-cover of $X$, there is $V_n\in\cV$ such that $F\cup P(\{0,\dots,\allowbreak m_n\})\sbst V_n$.
Take $U_n\in\cU$ such that $V_n^e\sbst U_n$. Then $U_n\nin\{U_0,\dots,U_{n-1}\}$.
As $V_n$ is open, for each $s\sbst \{0,\dots,m_n\}$ there is $k_s$ such that
for each $x\in P(\N)$ with $x\cap\{0,\dots,k_s-1\}=s$, $x\in V_n$. Let $m_{n+1}=\max\{k_s : s\sbst\{0,\dots,m_n\}\}$.

If $x_i\cap (m_n,m_{n+1})=\emptyset$ for all $i<e$, then $(x_0,\dots,x_{e-1})\in V_n^e\sbst U_n$.
\epf

The assumption in the theorem that there is an unbounded tower
of cardinality $\fb_0$ implies that $\fb_0=\fb$.
The proof is by induction on the power $e$ of $T\cup\Fin$.
The case $e=1$ follows from Theorem \ref{Sch}.

Let $\cU_0,\cU_1,\dots\in\Gamma((T\cup\Fin)^e)$. We may assume that these covers are countable.
As in the proof of Theorem \ref{Sch} (this time using Lemma \ref{GM++}),
there are for each $k$ $m^k_0<m^k_1<\dots$ and
 $U^k_0,U^k_1,\dots\in\cU_k$
(so that $\sseq{U^k_n}\in\Gamma((T\cup\Fin)^e)$), such that for
all $y_0,\dots,y_{e-1}\sbst\N$, $(y_0,\dots,y_{e-1})\in U^k_n$ whenever
$y_i\cap (m^k_n,m^k_{n+1})=\emptyset$ for all $i<e$.

Let $\alpha_0$ be such that $X_{\alpha_0}^e$ is not
contained in any member of $\Union_n\cU_n$.
As $T$ is unbounded, there is $\alpha$ such
that $\alpha_0\le \alpha<\fb$, and for each $k$,
$I_k=\{n : x_\alpha\cap (m^k_n,m^k_{n+1})=\emptyset\}$ is infinite.

Let $Y=\{x_\beta : \beta\ge\alpha\}$.
$(T\cup\Fin)^e\sm Y^e$ is a union of fewer than $\fb_0$ homeomorphic copies of $(T\cup\Fin)^{e-1}$.
By the induction hypothesis, $(T\cup\Fin)^{e-1}$ satisfies $\sgg$, and therefore so does
$(T\cup\Fin)^e\sm Y^e$.
For each $k$, $\{U^k_n : n\in I_k\}$ is a point-cofinite cover of $(T\cup\Fin)^e\sm Y^e$,
and thus there are infinite $J_0\sbst I_0, J_1\sbst I_1,\dots$, such that
$\{\bigcap_{n\in J_k}U^k_n : k\in\N\}$ is a point-cofinite cover of $(T\cup\Fin)^e\sm Y^e$.\footnote{Choosing
infinitely many elements from each cover, instead of one, can be done by adding to the given
sequence of covers all cofinite subsets of the given covers.}
For each $k$, pick $n_k\in J_k$ such that: $m^k_{n_k}>m^{k-1}_{n_{k-1}+1}$, $x_\alpha\cap (m^k_{n_k},m^k_{n_k+1})=\emptyset$,
and $U^k_{n_k}\nin\{U^0_{n_0},\dots,U^{k-1}_{n_{k-1}}\}$.

$\{U^k_{n_k} : k\in\N\}\in\Gamma(T\cup\Fin)$:
If $x\in (T\cup\Fin)^e\sm Y^e$, then $x\in U^k_{n_k}$ for all but finitely many $k$.
If $x=(x_{\beta_0},\dots,x_{\beta_{e-1}})\in Y$, then $\beta_0,\dots,\beta_{e-1}\ge\alpha$,
and thus $x_{\beta_0},\dots,x_{\beta_{e-1}}\as x_\alpha$.
For each large enough $k$, $m^k_{n_k}$ is large enough, so that
$x_{\beta_i}\cap (m^k_{n_k},m^k_{n_k+1})\sbst x_\alpha\cap (m^k_{n_k},m^k_{n_k+1})=\emptyset$
for all $i<e$, and thus $x\in U^k_{n_k}$.
\epf


There is an additional way to obtain nontrivial $\sgg$ sets:
The hypothesis $\fb=\cov(\cN)=\cof(\cN)$ provides $\fb$-Sierpi\'nski sets,
and $\fb$-Sierpi\'nski sets satisfy $\sgg$, even for \emph{Borel} point-cofinite
covers. Details are available in \cite{CBC}.

We record the following consequence of Theorem \ref{main} for later use.

\bcor\label{maincor}
For each unbounded tower $T$ of cardinality $\fb$, $T\cup\Fin$ satisfies $\sgg$ for open covers, but not
for Borel covers.
\ecor
\bpf
The latter property is hereditary for subsets \cite{CBC}.
By a theorem of Hurewicz, a set of reals satisfies $\ufin(\Gamma,\Gamma)$ if, and only if, each continuous
image of $X$ in $\NN$ is bounded. It follows that the set $T\sbst T\cup\Fin$ does not even satisfy $\ufin(\Gamma,\Gamma)$.
\epf

\section{A consistency result}\label{con}

By the results of the previous section, we have the following.

\blem
Assume that every set of
reals with property $\sgg$ has cardinality $<\fb$, and $\fc=\aleph_2$.
Then $\aleph_1=\ft=\cov(\cN)<\fb=\aleph_2$.
\elem
\bpf
As there is no unbounded tower, we have that $\ft<\fb=\fd$.
As $\fc=\aleph_2$, $\aleph_1=\ft<\fb=\aleph_2$. Since there are no $\fb$-Sirepi\'nski sets
and $\fb=\cof(\cN)=\fc$, $\cov(\cN)<\fb$.
\epf

In Laver's model \cite{Laver}, $\aleph_1=\ft=\cov(\cN)<\fb=\aleph_2$.
We will show that indeed, $\sgg$ is trivial there.
Laver's model was constructed to realize Borel's Conjecture,
asserting that ``strong measure zero'' is trivial.
In some sense, $\sgg$ is a dual of strong measure
zero. For example, the canonical examples of $\sgg$ sets are Sierpi\'nski sets, a measure
theoretic object, whereas the canonical examples of strong measure zero sets are Luzin sets, a Baire category
theoretic object. More about that can be seen in \cite{CBC}.

The main technical result of this paper is the following.

\bthm\label{bddLaver}
In the Laver model, if $X\sbst 2^\w$ has cardinality $\fb$,
then there is a Borel map $f:2^\w \to \w^\w$ such that $f[X]$ is unbounded.
\ethm
\bpf
The notation in this proof is as in Laver \cite{Laver}.
We will use the following slightly simplified version of Lemma 14 of \cite{Laver}.

\blem[Laver]\label{LaverLemma}
Let $\bbP_{\w_2}$ be the countable support iteration
of Laver forcing, $p\in\bbP_{\w_2}$, and $\name{a}$ be a
$\bbP_{\w_2}$-name such that
$$p\forces \name{a}\in 2^\w.$$
Then there are a condition $q$ stronger than $p$, and
finite $U_s\sbst 2^\w$ for each $s\in q(0)$ extending the root of $q(0)$, such that
for all such $s$ and all $n$:
$$q(0)_t\concat q\res [1,\w_2)\forces
\mbox{ `` }\exists u\in \check{U}_s \;\;u\res n =\name{a}
\res n \mbox{ '' }$$
for all but finitely many immediate successors $t$ of $s$
in $q(0)$.
\elem

Assume that $X\sbst 2^\w$ has no unbounded Borel image in
$\cM[G_{\w_2}]$, Laver's model.
For every code $u\in 2^\w$ for a Borel function
$f:2^\w\to\w^\w$
there exists $g\in\w^\w$ such that for every $x\in X$ we
have that $f(x)\leq^*g$.

By a standard L\"owenheim-Skolem argument, see Theorem 4.5 on page 281 of \cite{baum},
or section 4 on page 580 of \cite{onto}, we may find $\alpha<\w_2$ such that for every
code $u\in \cM[G_\alpha]$ there is an upper bound $g\in\cM[G_\alpha]$.
By the arguments employed by
Laver \cite[Lemmata 10 and 11]{Laver}, we may assume that $\cM[G_\alpha]$ is the ground model $\cM$.

Since the continuum hypothesis holds in $\cM$ and $|X|=\fb=\aleph_2$,
there are $p\in G_{\w_2}$ and $\name{a}$ such that
$$p\forces \name{a}\in \name{X} \rmand \name{a}\notin \cM.$$
Work in the ground model $\cM$.

Let $q\leq p$ be as in Lemma \ref{LaverLemma}.
Define
$$Q=\{s\in q(0)\;:\; \op{root}(q(0))\sbst s\}$$ and
let $U_s$, $s\in Q$, be the finite sets from the Lemma.
Let $U=\bigcup_{s\in Q}U_s$.
Define a Borel map $f:2^\w\to \w^Q$ so that
for every $x\in 2^\w\sm U$ and
for each $s\in Q$: If $f(x)(s)=n$, then $x\res n\neq u\res n$
for each $u\in U_s$. For $x\in U$, $f(x)$ may be arbitrary.
There must be a $g\in \w^Q\cap \cM$ and $r\leq q$ such that
$$r\forces f(\name{a})\leq^* \check{g}.$$
Since $p$ forced that $a$ is not in the ground model, it
cannot be that $a$ is in $U$.  We may extend $r(0)$ if necessary
so that if $s=\op{root}(r(0))$, then
$$r\forces f(\name{a})(s)\leq \check{g}(s).$$
But this is a contradiction to Lemma \ref{LaverLemma}, since for all but finitely
many $t\in r(0)$ which are immediate extensions of $s$:
$$r(0)_t\concat q\res [1,\w_2)\forces f(\name{a})(s)>\check{g}(s).$$
\qedhere
\epf

In \cite{hH}, Tsaban and Zdomskyy prove that $\sgg$ for Borel covers is equivalent to
the Ko\v{c}inac property $\scof(\Gamma,\Gamma)$ \cite{KocAlpha},
asserting that for all $\cU_0,\cU_1,\dots\in\Gamma$, there are
cofinite subsets $\cV_0\sbst\cU_0,\cV_1\sbst\cU_1,\dots$ such that
 $\Union_n\cV_n\in\Gamma$.
The main result of \cite{Dow90} can be reformulated as follows.

\bthm[Dow \cite{Dow90}]\label{dowthm}
In Laver's model, $\sgg$ implies $\scof(\Gamma,\Gamma)$.
\ethm
For the reader's convenience, we give Dow's proof,
adapted to the present notation.
\bpf
A family $\cH\sbst\roth$ is \emph{$\om$-splitting} if for each countable $\cA\sbst\roth$,
there is $H\in\cH$ which splits each element of $\cA$, i.e.,
$$|A\cap H|= |A\sm H|=\om \mbox{ for all } A\in\cA.$$
The main technical result in \cite{Dow90} is the following.

\blem[Dow]\label{hitlem}
In Laver's model, each $\om$-splitting family contains an
$\om$-splitting family of cardinality $<\fb$.
\elem

\noindent Assume that $X$ satisfies $\sgg$.
Let $\cU_0,\cU_1,\dots$
be open point-cofinite countable covers of $X$.
We may assume\footnote{To see why, replace each $\cU_n$ by
$\cU_n\sm\bigcup_{i<n}\cU_i$, and discard the finite ones.
It suffices to show that $\scof(\Gamma,\Gamma)$ applies
to those that are left.}
that $\cU_i\cap\cU_j=\emptyset$ whenever $i\neq j$.
Put $\cU=\bigcup_{n<\om}\cU_n$.  We identify
$\cU$ with $\om$, its cardinality.

Define $\cH\su [\cU]^\om$ as follows.
For $H\in [\cU]^\om$, put $H\in\cH$ if and only if there exists $\cV\in[\cU]^\om$,
a point-cofinite cover of $X$, such that $H\cap \cU_n\su^* \cV$ for all $n$.
We claim that $\cH$ is an $\om$-splitting family.
As $\cH$ is closed under taking infinite subsets, it suffices to show that
it is \emph{$\om$-hitting}, i.e., for any countable $\cA\su [\cU]^\om$
there exists $H\in \cH$ which intersects each $A\in\cA$.
(It is enough to intersect each $A\in\cA$, since we may assume that
$\cA$ is closed under taking cofinite subsets.)

Let $\cA\su [\cU]^\om$ be countable.
For each $n$, choose sets $\cU_{n,m}\in [\cU_n]^\om$, $m\in\om$,
such that for each $A\in\cA$, if $A\cap \cU_n$ is
infinite, then $\cU_{n,m}\su A$ for some $m$.
Apply the $\sgg$ to the family $\{\cU_{n,m}:n,m\in\N\}$,
to obtain a point-cofinite $\cV\su\cU$ such that
$\cV\cap \cU_{n,m}$ is nonempty for all $n,m$.

Next, choose finite subsets $\cF_n\su \cU_n$, $n\in\N$,
such that for each $A\in \cA$ with $A\cap \cU_n$ finite for all $n$,
then $A\as\Union_n\cF_n$.
Take $H=\cV\cup\bigcup_{n}\cF_n$.
Then $H$ is in $\cH$ and meets each $A\in \cA$.
This shows that $\cH$ is an $\om$-splitting family.

By Lemma \ref{hitlem}, there is an $\om$-splitting $\cH'\sbst\cH$ of
cardinality $<\fb$. For each $H\in\cH'$, let $\cV_H$ witness that
$H$ is in $\cH$, i.e., $\cV_H\su\cU$ is a point-cofinite cover of $X$ and
$H\cap \cU_n\su^* \cV_H$ for all $n$.

By the definition of $\fb$, we may find finite $\cF_n\su\cU_n$, $n\in\om$,
such that for each $H\in \cH'$,
$$H\cap \cU_n\su \cV_H\cup\cF_n$$
for all but finitely many $n$.
We claim that  $\cW=\bigcup_{n}\cU_n\sm\cF_n$ is point-cofinite.
Suppose it is not. Then there is $x\in X$ such that for infinitely
many $n$, there is $U_n\in \cU_n\sm\cF_n$ with
$x\notin U_{n}$.  Let $H\in \cH'$ contain infinitely
many of these $U_n$.
By the above inclusion, all but finitely many
of these $U_n$ are in $\cV_H$.  This contradicts the fact
that $\cV_H$ is point-cofinite.
\epf

We therefore have the following.

\bthm\label{mainthm}
In Laver's model, each set of reals $X$ satisfying $\sgg$ has cardinality less
than $\fb$.
\ethm
\bpf
By Dow's Theorem, $\sgg$ implies $\scof(\Gamma,\Gamma)$, which in turn implies $\sgg$ for Borel covers \cite{hH}.
The latter property is equivalent to having all Borel images in $\NN$ bounded \cite{CBC}. Apply Theorem \ref{bddLaver}.
\epf

Thus, it is consistent that strong measure zero and $\sgg$ are both trivial.

\medskip

The proof of Dow's Theorem \ref{dowthm} becomes more natural after replacing,
in Lemma \ref{hitlem} ``$\w$-splitting'' by ``$\w$-hitting''.
This is possible, due to the following fact (cf.\ Remark 4 of \cite{Dow90}).

\bprp
For each infinite cardinal $\kappa$, the following are equivalent:
\be
\itm Each $\om$-splitting family contains an $\om$-splitting family of cardinality $<\kappa$.
\itm Each $\om$-hitting family contains an $\om$-hitting family of cardinality $<\kappa$.
\ee
\eprp
\bpf
$(1\Impl 2)$
Suppose $\cA$ is an $\om$-hitting family.  Let
$\cB=\bigcup_{A\in\cA}[A]^\om$.
Then $\cB$ is $\om$-splitting.  By (1) there exists
$\cC\sbst\cB$ of size $<\kappa$ which is $\om$-splitting.
Choose $\cD\sbst\cA$ of size $<\kappa$ such that for every
$C\in\cC$ there exists $D\in\cD$ with $C\sbst D$.
Then $\cD$ is $\om$-hitting.

$(2\Impl 1)$ Suppose $\cA$ is an $\om$-splitting family.  For each
$A\sbst\om$ define
$$A^*=\{2n:n\in A\}\cup\{2n+1:n\in \overline{A}\}.$$
Then the family $\cA^*=\{A^*:A\in\cA\}$ is $\om$-hitting. To see this, suppose that $\cB$ is countable.
Without loss we may assume that $\cB=\cB_0\cup \cB_1$
where each element of $\cB_0$ is
a subset of the evens and each element of $\cB_1$
is a subset of the odds.
For $B\in\cB_0$ let
$C_B=\{n:2n\in B\}$ and for $B\in\cB_1$ let
$C_B=\{n:2n+1\in B\}$.
Now put
$$\cC=\{C_B:B\in\cB\}.$$
Since $\cA$ is $\om$-splitting
there is $A\in\cA$ which splits $\cC$.  If
$B\in B_0$, then $A\cap C_B$ infinite implies
$B\cap A^*$ infinite.  If $B\in B_1$ then
$\overline{A}\cap C_B$ infinite implies $B\cap A^*$ infinite.

By (2) there exists $\cA_0\sbst \cA$ of cardinality $<\kappa$
such that $\cA_0^*$ is $\om$-hitting.  We claim that
$\cA_0$ is $\om$-splitting.  Given any
$B\sbst\om$ let $B^\prime=\{2n:n\in B\}$
and let $B^{\prime\prime}= \{2n+1:n\in B\}$.  Given
$\cB\sbst[\om]^\om$ countable,
there exists $A\in\cA_0$ such
that $A^*$ hits each $B^\prime$ and $B^{\prime\prime}$ for $B\in\cB$.
But this implies that $A$ splits $B$.
\epf

\section{Applications to Arhangel'ski\u{\i}'s $\alpha_i$ spaces}

Let $Y$ be a general (not necessarily metrizable) topological space.
We say that a countably infinite set $A\sbst Y$ converges to a point $y\in Y$ if
each (equivalently, some) bijective enumeration of $A$ converges to $y$.
The following concepts are due to Arhangel'ski\u{\i}.
$Y$ is an \emph{$\alpha_1$ space} if for each $y\in Y$ and each
sequence $A_0,A_1,\dots$ of countably infinite sets, each converging to $y$,
there are cofinite $B_0\sbst A_0, B_1\sbst A_1,\dots$, such that $\Union_nB_n$ converges
to $y$. Replacing ``cofinite'' by ``singletons'' (or equivalently, by ``infinite''),
we obtain the definition of an \emph{$\alpha_2$ space}.

We first consider countable spaces.

\bdfn\label{local}
Let $X$ be a set of reals, and let $\cU_0,\cU_1,\dots$ be countable point-cofinite covers of $X$.
For each $n$, enumerate bijectively $\cU_n=\{U^n_m : m\in\N\}$.
We associate to $X$ a (new) topology $\tau$ on the fan $S_\w=\N\x\N\cup\{\oo\}$
as follows: $\oo$ is the only nonisolated point of $S_\w$,
and a neighborhood base at $\oo$ is given by the sets
$$[\oo]_F = \{(n,m) : F\sbst U^n_m\}$$
for each finite $F\sbst X$.
\edfn

\blem\label{dict}
In the notation of Definition \ref{local}:
$A$ converges to $\oo$ in $\tau$ if, and only if,
$\cU(A)=\{U^n_m : (n,m)\in A\}$ is a point-cofinite cover of $X$.\qed
\elem

Assume that there is an unbounded tower. By Corollary \ref{maincor},
there is a set of reals $X$ satisfying $\sgg$ but not $\scof(\Gamma,\Gamma)$.
Let $\cU_0,\cU_1,\dots$ be countable open point-cofinite covers of $X$ witnessing
the failure of $\scof(\Gamma,\Gamma)$.
Then, by Lemma \ref{dict}, $(S_\w,\tau)$ is $\alpha_2$ but not $\alpha_1$.
In particular, we reproduce the following.

\bcor[Nyikos \cite{Nyikos92}]\label{Nyik}
If there is an unbounded tower of cardinality $\fb$,
then there is a countable $\alpha_2$ space, which is not an $\alpha_1$ space.\qed
\ecor

Recall that by Proposition \ref{anytower}, it suffices to assume in Corollary \ref{Nyik} the existence
of any unbounded tower.

\medskip

Next, we consider spaces of continuous functions.
Consider $C(X)$, the family of continuous real-valued functions, as a subspace of the Tychonoff
product $\bbR^X$, i.e., with the topology of pointwise convergence.
Sakai \cite{Sakai07} proved that $X$ satisfies $\sgg$
for \emph{clopen} covers if,
and only if, $C(X)$ is an $\alpha_2$ space.
The main result of \cite{hH} is that $C(X)$ is $\alpha_1$ if, and only if,
$X$ satisfies $\sgg$ for Borel covers (equivalently, each Borel image of $X$ in $\NN$ is bounded).

The \emph{Scheepers Conjecture} is that for subsets of $\bbR\sm\bbQ$,
$\sgg$ for clopen covers implies $\sgg$ for open covers.
Dow \cite{Dow90} proved that in Laver's model, every $\alpha_2$ space is $\alpha_1$.
By Theorem \ref{bddLaver}, we can add the last item in the following list.

\bcor\label{morecor}
In Laver's model, the following are equivalent for sets of reals $X$:
\be
\itm $C(X)$ is an $\alpha_2$ space;
\itm $C(X)$ is an $\alpha_1$ space;
\itm $X$ satisfies $\sgg$ for clopen covers;
\itm $X$ satisfies $\sgg$ for open covers;
\itm $X$ satisfies $\sgg$ for Borel covers;
\itm $|X|<\fb$.\qed
\ee
\ecor

On the other hand, Corollary \ref{maincor} implies the following.

\bcor\label{last}
If there is an unbounded tower, then there is a set of reals $X$ such that $C(X)$ is $\alpha_2$ but not $\alpha_1$.\qed
\ecor

Essentially, Corollary \ref{Nyik} is a special case of Corollary \ref{maincor},
whereas Corollary \ref{last} is equivalent to Corollary \ref{maincor}.

\subsection*{Acknowledgment}
We thank Lyubomyr Zdomskyy for his useful comments.

\ed
\begin{thebibliography}{99}
\bibitem{BaSh01}
T.\ Bartoszy\'nski and S.\ Shelah,
\emph{Continuous images of sets of reals},
Topology and its Applications \textbf{116} (2001), 243--253.

\Pa{ideals}{T. Bartoszy\'nski and B. Tsaban}{Hereditary topological diagonalizations and the Menger-Hurewicz Conjectures}{Proceedings of the American Mathematical Society}{134}{2006}{605}{615}

\bibitem{baum}
J. Baumgartner, R. Laver,
\emph{Iterated perfect-set forcing},
Annals of Mathematical Logic \textbf{17} (1979), 271--288.

\bibitem{BlassHBK}
A. Blass,
\emph{Combinatorial cardinal characteristics of the continuum},
in: \textbf{Handbook of Set Theory} (M.\ Foreman, A.\ Kanamori, and M.\ Magidor, eds.),
Kluwer Academic Publishers, Dordrecht, to appear.
\texttt{http://www.math.lsa.umich.edu/\~{}ablass/hbk.pdf}

\bibitem{Dow90}
A. Dow,
\emph{Two classes of Fr\'echet-Urysohn spaces},
Proceedings of the American Mathematical Society \textbf{108} (1990), 241--247.

\bibitem{GM}
F.\ Galvin and A. Miller,
\emph{$\gamma$-sets and other singular sets of real numbers},
Topology and its Applications \textbf{17} (1984), 145--155.

\bibitem{GN}
J.\ Gerlits and Zs.\ Nagy,
\emph{Some properties of $C(X)$, I},
Topology and its Applications \textbf{14} (1982), 151--161.

\bibitem{Hure25}
W.\ Hurewicz,
\emph{\"Uber eine Verallgemeinerung des Borelschen Theorems},
Mathematische Zeitschrift \textbf{24} (1925), 401--421.

\bibitem{coc2}
W. Just, A. Miller, M. Scheepers, and P. Szeptycki,
\emph{The combinatorics of open covers II},
Topology and its Applications \textbf{73} (1996), 241--266.

\bibitem{KocAlpha}
L.\ Ko\v{c}inac,
\emph{Selection principles related to $\alpha_i$-properties},
Taiwanese Journal of Mathematics 12 (2008), 561--572.

\bibitem{Laver}
R.\ Laver,
\emph{On the consistency of Borel's conjecture},
Acta Mathematica \textbf{137} (1976), 151--169.

\bibitem{onto}
A. Miller,
\emph{Mapping a set of reals onto the reals},
Journal of Symbolic Logic \textbf{48} (1983), 575--584.

\bibitem{Nyikos92}
P. Nyikos,
\emph{Subsets of $\omega^\omega$ and the Fr\'echet-Urysohn and $\alpha_i$-properties},
Topology and its Applications \textbf{48} (1992), 91--116.

\bibitem{Sakai07}
M.\ Sakai,
\emph{The sequence selection properties of $C_p(X)$},
Topology and its Applications \textbf{154} (2007), 552--560.

\bibitem{coc1}
M.\ Scheepers,
\emph{Combinatorics of open covers I: Ramsey theory},
Topology and its Applications \textbf{69} (1996), 31--62.

\bibitem{alpha_i}
M.\ Scheepers,
\emph{$C_p(X)$ and Arhangel'ski\u{\i}'s $\alpha_i$ spaces},
Topology and its Applications \textbf{89} (1998), 265--275.

\bibitem{wqn}
M.\ Scheepers,
\emph{Sequential convergence in ${\sf C}_p(X)$ and a covering property},
East-West Journal of Mathematics \textbf{1} (1999),
207--214.

\Pa{CBC}{M. Scheepers and B. Tsaban}{The combinatorics of Borel covers}{Topology and its Applications}{121}{2002}{357}{382}

\bibitem{MHP}
B. Tsaban,
\emph{Menger's and Hurewicz's Problems: Solutions from ``The Book'' and refinements},
Contemporary Mathematics, to appear.


\bibitem{hH}
B. Tsaban and L. Zdomskyy,
\emph{Hereditarily Hurewicz spaces and Arhangel'ski\u{\i} sheaf amalgamations},
Journal of the European Mathematical Society, to appear.

\end{thebibliography}
